\documentclass[11pt]{amsart}

\usepackage{amsmath}
\usepackage{amssymb}
\usepackage{mathrsfs}
\usepackage[all]{xy}

\newtheorem{theorem}{Theorem}[section]
\newtheorem{claim}[theorem]{Claim}

\theoremstyle{definition}
\newtheorem{definition}[theorem]{Definition}

\theoremstyle{remark}

\newcount\skewfactor
\def\mathunderaccent#1#2 {\let\theaccent#1\skewfactor#2
\mathpalette\putaccentunder}
\def\putaccentunder#1#2{\oalign{$#1#2$\crcr\hidewidth
\vbox to.2ex{\hbox{$#1\skew\skewfactor\theaccent{}$}\vss}\hidewidth}}
\def\name{\mathunderaccent\tilde-3 }

\def\smallbox#1{\leavevmode\thinspace\hbox{\vrule\vtop{\vbox
   {\hrule\kern1pt\hbox{\vphantom{\tt/}\thinspace{\tt#1}\thinspace}}
   \kern1pt\hrule}\vrule}\thinspace}

\def\qedref#1{$\qed_{\reforiginal{#1}}$}

\title{Yablo's paradox and forcing}
\author{Shimon Garti}
\address{Einstein Institute of Mathematics
 The Hebrew University of Jerusalem,
 Jerusalem 91904, Israel}
\email{shimon.garty@mail.huji.ac.il}


\begin{document}
\let\labeloriginal\label
\let\reforiginal\ref
\def\ref#1{\reforiginal{#1}}
\def\label#1{\labeloriginal{#1}}

\begin{abstract}
We discuss the problem of self-reference in Yablo's paradox from the point of view of the relationship between names and objects.
For this end, we introduce a forcing version of the paradox and try to understand its implication on the self-referential component of the paradox.
\end{abstract}

\maketitle

\newpage

\section{Introduction}

Yablo introduced in \cite{MR1249561} a paradox which allegedly belongs to the liar paradox family yet contains no self-reference.
The paradox consists of an infinite family of sentences $\mathcal{S}=\{S_n:n\in\omega\}$ where each $S_n$ is the statement $\forall\ell>n,S_\ell$ is untrue.
A moment of perusal leads to a contradiction which forms the paradox.
A central theme concerning Yablo's paradox is the correctness of the title of Yablo's paper.
According to the title, this paradox is not self-referential.
Whether this is true or not is still under debate, see for example \cite{MR1837271}, \cite{MR1607031} and \cite{MR2382232}.

The main claim of the present paper is that this question depends on another important philosophical issue, the relationship between mathematical objects and their names.
We will try to examine Yablo's paradox from the point of view of forcing theory.
Being a theory which enables us to separate objects from their names in an accurate way, we believe that it may shed light on the paradox in general and on the self-referential component connected with it.
We shall use Prikry forcing, which appeared in \cite{MR262075}.
We suggest \cite{MR2768695} as a source for background on Prikry forcing and we employ the Jerusalem forcing notation as done in this monograph.

\section{Some remarks on forcing}

Paul Cohen introduced the method of forcing in \cite{MR0232676}, in order to prove the consistency of the failure of the continuum hypothesis and the axiom of choice.
The modern practice of forcing owes a lot to Shoenfield, \cite{MR0280359}.
In this section we describe, shortly, the basic idea of this method.

We begin with a universe of set theory which is called \emph{the ground model} and denoted by $V$.
We define, in $V$, a partial order $\mathbb{P}$.
The set $\mathbb{P}$ is called a forcing notion and the elements of $\mathbb{P}$ are called conditions.
Each condition gives a small piece of information about a mathematical object that we try to force.

This object does not belong to $V$, but its possible existence does not contradict the axioms of set theory.
From $V$ we can imagine the desired object, give it a name and describe it to some extent.
In the next step we choose a generic set $G\subseteq\mathbb{P}$, which does not belong to $V$ (unless $\mathbb{P}$ is trivial in some sense).
Now we extend the universe $V$ by adding $G$ and also many other sets.

The resulting model is called \emph{the generic extension}, and denoted by $V[G]$.
There are two ways to describe the process of enlarging $V$ to $V[G]$.
The axiomatic description says that we add $G$ and then close the universe under the axioms, so we must add subsets of $G$, applications of the replacement axiom in which $G$ is involved, and so on.
An alternative and more practical way to describe $V[G]$ is through the concept of \emph{names}.
Every object $S\in V[G]$ has a name $\name{S}$ in the ground model $V$.
A recursive definition provides the ability to interpret names using the generic set $G$.
Thus $V[G]$ is simply the interpretation of all the $\mathbb{P}$-names from the ground model according to the generic set $G$.

This process is somewhat parallel to the geometric approach of the ancient Greeks.
Let $\mathcal{F}$ be the field of numbers constructible by a ruler and a pair of compasses, the ground model of the Greeks.
We know that $\pi\notin\mathcal{F}$, though the Greeks did not know this.
Nonetheless, they could name it and describe it quite vividly as the area of a circle whose radius is 1.
Moreover, they could give partial information about $\pi$ using areas of bounded (and bounding) polygons.
The most important thing is that an extension of their universe embodies $\pi$ as a real number.

This is quite similar to the forcing process.
One begins with a ground model and a desired mathematical object.
This object is not an element of the ground model, but it can be described from the ground model.
In the official terminology, it has a name in the ground model $V$.
By extending $V$ this name can be interpreted and become a real object in the generic extension $V[G]$.

The object that we would like to force is the set of sentences $\mathcal{S}$ which forms Yablo's paradox.
In the ground model we will define conditions (and a partial order between them), each of which gives only partial information about $\mathcal{S}$.
This construction will be depicted in the next section.

\section{Prikry forcing and Yablo's paradox}

Consider the statement $S_n$ which says $\forall\ell>n,\neg S_\ell$.
Although $S_n$ does not refer to $S_n$, it mentions the index $n$.
The explicit reference to the index seems inevitable in any formulation of Yablo's paradox.
We are asking, therefore, whether such a reference is a self-reference.

From this point of view, the question of self-reference in Yablo's paradox is a reflection of the substantial issue of names and objects.
If one argues that Yablo's paradox is not self-referential then one claims that the \emph{name} of an object differs significantly from the object itself.
Therefore, a reference to $n$, the name of $S_n$, is not necessarily a self-reference.
If one identifies names and objects by claiming that a (mathematical) object is nothing but the union of its names, then one concludes that Yablo's paradox is self-referential indeed.
Observe that this point also appears in other formulations of Yablo's paradox, like the set-theoretical version of Goldstein in \cite{MR1324809}.

We suggest below a formal angle from which this aspect of self-reference in Yablo's paradox can be examined.
This will be done, basically, by an application of forcing.
We shall force the existence of the set $\mathcal{S}$, and we wish to say something about self-reference with respect to this set.

It has been claimed that though each $S_n$ is not self-referential, the entire collection $\mathcal{S}$ contains a self-reference.
In order to build $\mathcal{S}$ one has to define every $S_n$.
In order to define $S_n$ one has to say something about an end-segment of $\mathcal{S}$.
Thus the existence of $\mathcal{S}$ captures some circularity, and here lies the quintessential self-reference of the paradox, see for example the short observation of Smith in \cite{MRreview}.
We indicate that this point is also a central issue when one tries to understand why Yablo's paradox does not lead to a contradiction in formal systems of set theory like \textsf{ZFC}.

Our goal is to phrase a version of Yablo's paradox in which the statement $S_n$ does not refer to $\mathcal{S}$, or to an end-segment of $\mathcal{S}$.
In our formulation, each $S_n$ mentions only a finite set of $S_\ell$s.
To gather all the information we will extend our universe of set theory.

\begin{definition}
\label{defforcing} Yablo's paradox forcing notion. \newline 
Let $\mathscr{U}$ be a normal ultrafilter over a measurable cardinal $\kappa$.
We define a Prikry-type forcing notion $\mathbb{P}$.
\begin{enumerate}
\item [$(\aleph)$] A condition $p\in\mathbb{P}$ is a triple $(s,S,A)=(s^p,S^p,A^p)$ where $s\in[\kappa]^{<\omega}, A\in\mathscr{U}$ and $\max(s)<\min{A}$. If $s=\{\alpha_0,\ldots,\alpha_{n-1}\}$ and $i<j\Rightarrow\alpha_i<\alpha_j$ then $S=\{S_{\alpha_0},\ldots,S_{\alpha_{n-1}}\}$ where $S_{\alpha_i}=\bigwedge_{i<j\leq n}\neg S_{\alpha_j}$.
\item [$(\beth)$] If $p,q\in\mathbb{P}$ then $p\leq q$ iff $s^p\trianglelefteq s^q$ and $A^q\subseteq A^p$.
\item [$(\gimel)$] If $p,q\in\mathbb{P}$ then $p\leq^*q$ iff $s^p=s^q, S^p=S^q$ and $A^q\subseteq A^p$.
\end{enumerate}
\end{definition}

The forcing notion $\mathbb{P}$ is a variant of Prikry forcing, augmented with the elements of the form $S^p$ which approximate Yablo's system as will be depicted anon.
Ahead of this description let us indicate that this forcing is a Prikry-type forcing notion.

\begin{claim}
\label{clmproperties} The forcing notion $(\mathbb{P},\leq,\leq^*)$ is $\kappa^+$-cc, satisfies Prikry property and $(\mathbb{P},\leq^*)$ is $\kappa$-closed.
Hence all cardinals are preserved in the generic extension by $\mathbb{P}$.
\end{claim}

\par\noindent\emph{Proof}. \newline 
We commence with the chain condition.
If $p,q\in\mathbb{P}$ and $s^p=s^q$ then necessarily $S^p=S^q$ and it follows that $p\parallel q$, so the size of any antichain is bounded by the number of stems which is $\kappa$.
Prikry property can be proved using Rowbottom's theorem as done with the usual Prikry forcing.
The closure degree of $(\mathbb{P},\leq^*)$ is $\kappa$ since $\mathscr{U}$ is normal and in particular $\kappa$-complete.
Combining all the above statements we see that all cardinals are preserved, as required.

\hfill \qedref{clmproperties}

Let $G\subseteq\mathbb{P}$ be generic over $V$ and let $(\rho_n:n\in\omega)$ be the associated Prikry sequence added by the generic set.
Define $\mathcal{S}=\bigcup\{S^p:p\in G\}$.
One can verify that $\mathcal{S}$ is a system of $\omega$-many sentences $\{S_{\rho_n}:n\in\omega\}$ which exemplifies Yablo's paradox in $V[G]$.

Let us add a few words about the connection between $\mathcal{S}\in V[G]$ and the conditions of $\mathbb{P}$.
Every $p\in\mathbb{P}$ contains a finite approximation $S^p$ of $\mathcal{S}$.
For each $p$, the set $S^p$ has a last element and hence $S^p$ forms no paradox.
If $q$ is another condition then possibly $p\perp q$, so $q$ contains a different (mayhap contradictory) information about $\mathcal{S}$.
However, $G$ is directed and hence picks only compatible conditions for creating $\mathcal{S}$ in the generic extension.

The set $\mathcal{S}$ is an object in $V[G]$, but not in the ground model $V$.
In the next section we shall discuss the meaning of this fact with respect to self-reference in Yablo's paradox.
For the present section we just indicate that $\mathcal{S}$ is not an object but it has a name $\name{\mathcal{S}}$ in $V$, hence the above \emph{generic} version of Yablo's paradox highlights the name-object aspect of the self-reference issue.

\section{The ground model and the generic extension}

Our main claim in this section is that the forcing version of Yablo's paradox is free from self-reference even if one argues that the original paradox is self-referential.
Remark that $\mathcal{S}\in V[G]$ and therefore if we study the paradox in $V[G]$ then it is the usual Yablo's paradox.
In particular, the debate and the various arguments concerning circularity and self-reference are in the same position.

But let us try to consider the paradox from the ground model's point of view.
It seems that in $V$ there is no self-reference at all, as each $S_n$ refers only to finitely many statements, all of which appear explicitly in the condition $p$ to which $S_n$ belongs.
In particular, there is no appeal to $\mathcal{S}$ or some end-segment of $\mathcal{S}$.
One may argue, however, that from the ground model's point of view there is no paradox at all.
As far as the set $\mathcal{S}$ does not exist we cannot obtain the desired contradiction, and this set does not exist in $V$.
One may agree that we have a name $\name{\mathcal{S}}$ for this set and hence we have \emph{a name of a paradox} but not an actual paradox.

Here we touch upon the issue of the definition of a paradox.
In its basic form, a paradox is a collection of statements which lead to a contradiction.
The forcing version introduced in the previous section certainly leads to a contradiction, and hence should be regarded as a paradox.
In order to get a contradiction one has to invoke some set-theoretical arguments, including the mathematical principles which establish the generic extension by $\mathbb{P}$.
Yet there is no essential difference between these principles and the parallel process rendered in deriving a contradiction at any other paradox.
Though we have only a name $\name{\mathcal{S}}$ in $V$, we know how to describe a formal process in which we choose a generic object $G$, interpret all the names including $\name{\mathcal{S}}$, and reach to a contradiction.
We conclude, therefore, that even in $V$ we have a paradox, this time with no self-reference.

\section{Acknowledgement}

I thank the referee of the paper for his/her work, especially for suggesting the incorporation of the second section into the paper.

\newpage 

\bibliographystyle{amsplain}
\bibliography{arlist}

\end{document}